\documentclass{article}

\usepackage{amsmath,amsfonts,amssymb,amsthm,cleveref,color}
\usepackage{enumerate}
\usepackage{ifthen}

\newtheorem{theorem}{Theorem}

\newtheorem{proposition}[theorem]{Proposition}

\newtheorem{observation}[theorem]{Observation}
\newtheorem{observations}[theorem]{Observations}
\newtheorem{remark}[theorem]{Remark}

\newtheorem{lemma}[theorem]{Lemma}
\newtheorem{problem}[theorem]{Problem}
\newtheorem{facts}[theorem]{Facts}

\begin{document}
\title{Effective MC-finiteness}
\author{Yuval Filmus, Eldar Fischer, Johann A. Makowsky \\ \small Faculty of Computer Science \\ \small Technion-Israel Institute of Technology, Haifa, Israel}

\maketitle
\begin{abstract}
An integer sequence $(a_n)_{n \in \mathbb{N}}$ is \emph{MC-finite} if for all $m$, the sequence $a_n \bmod m$ is eventually periodic.
There are MC-finite sequences $(a_n)_{n \in \mathbb{N}}$ such that the function $F: (m,n) \mapsto a_n \bmod m$ is not computable.
In \cite{filmus2023mc} 
%@article{filmus2023mc,
%  title={MC-Finiteness of Restricted Set Partition Functions},
%  author={Filmus, Yuval and Fischer, Eldar and Makowsky, Johann A and Rakita, Vsevolod},
%  journal={Journal of Integer Sequences},
%  volume={26},
%  number={2},
%  pages={3},
%  year={2023}
%}
we presented concrete examples of MC-finite sequences taken from the Online Encyclopedia of Integer Sequences (OEIS)
without discussing the  computability of $F$.
In this paper we discuss cases when this $F$ is  effectively computable. 

\end{abstract}

%======================================================================
\newcommand{\cC}{\mathcal{C}}
\newcommand{\cD}{\mathcal{D}}
\newcommand{\cE}{\mathcal{E}}
\newcommand{\cT}{\mathcal{T}}
\newcommand{\fA}{{\mathfrak A}}
\newcommand{\fB}{{\mathfrak B}}
\newcommand{\fC}{{\mathfrak C}}
\newcommand{\N}{{\mathbb N}}
\newcommand{\Z}{{\mathbb Z}}
\newcommand{\MSOL}{\mathbf{MSOL}}
\newcommand{\SOL}{\mathbf{SOL}}
\newcommand{\CMSOL}{\mathbf{CMSOL}}
\newcommand{\FOL}{\mathbf{FOL}}

%======================================================================
%\marginpar{File: FFM-draft-1}
%\input{FFM-draft-1}
%\marginpar{File: FFM-intro}
\section{Introduction}
\label{sec:introduction}
An integer sequence $(a_n)_{n \in \mathbb{N}}$ is \emph{MC-finite} if for all $m$, the sequence $a_n \bmod m$ is eventually periodic.

If an integer sequence is MC-finite, then for each $m$, we can compute the function $n \mapsto a_n \bmod m$ with an algorithm $A_m(n)$ 
in polynomial time, as follows. 
Suppose that $a_{n+r} \equiv a_n \pmod{m}$ for all $n \geq N$. We store in a look-up table the entries
\[
 a_1,\dots,a_{N+r-1}.
\]
The algorithm is as follows:
\begin{enumerate}
    \item If $n < N$, then output $a_n$.
    \item Otherwise, let $i = (n - N) \bmod m$. \\ Output $a_{N + i}$.
\end{enumerate}

In this note we will be interested in the computability of $a_n \bmod m$ as a function of both $m$ and $n$. In general, 
this function need not be computable, as the following example shows.

\begin{theorem} \label{thm:not computable}
There is an MC-finite sequence $(a_n)_{n \in \mathbb{N}}$ such that the function $(m,n) \mapsto a_n \bmod m$ is not-computable.
\end{theorem}
\begin{proof}
Let $(b_n)_{n \in \mathbb{N}}$ be any not computable $0,1$ sequence, and let $(p_n)_{n \in \mathbb{N}}$ be the sequence of primes. 
Using the Chinese remainder theorem, let $a_n$ be the unique integer in $\{1, \dots, p_1^n \cdots p_n^n\}$ 
such that $a_n \equiv b_k \pmod{p_k^n}$ for all $k \leq n$.

We claim that the sequence $a_n$ is MC-finite. 
By the Chinese remainder theorem, it suffices to show that it is eventually periodic modulo any prime power $p_i^j$. 
Indeed, for $n \geq \max(i,j)$ we have $a_n \bmod p_i^j = b_n$.

On the other hand, since $a_n \bmod p_n = b_i$, the function $(m,n) \mapsto a_n \bmod m$ is not computable.
\end{proof}

If $(a_n)_{n \in \mathbb{N}}$ is computable, then clearly $(m,n) \mapsto a_n \bmod m$ is also computable, 
and so the interesting question is the complexity of computing $a_n \bmod m$ as a function of $m$ and $n$. 

Assume you look at an integer sequence $a_n$ of which the values of $a_n: n \leq n_0$ are listed but
\begin{itemize}
\item
you do not have the resources to recompute these values, and
\item
you suspect that for $k \leq n_0$ the value $a_k =b$ is wrong.
\end{itemize}
If you know that $a_n$ is MC-finite and you can establish that $a_k \not \equiv b \bmod m$ this would confirm your suspicion.
In other words, an algorithm which computes $(m,n) \mapsto a_n \bmod m$ is a {\em sufficient falsification method} for checking
$a_n =b$ for MC-finite seuqences $a_n$.

Indeed, the concept (but not its name) of MC-finiteness was introduced in 
\cite{pr:BlatterSpecker84,specker1990application} as such a falsification method. A systematic study of MC-finiteness is given in \cite{filmus2023mc},
however, without discussing complexity issues.

We say that the problem is \emph{fixed-parameter tractable} (FPT) if $a_n \bmod m$ can be computed in time $T(m) \cdot n^{O(1)}$, 
where $T(m)$ is an arbitrary function.
We say that a sequence $a_n$ is effectively MC-finite if computing $a_n \bmod m$ is FPT.

%--------------------------------------
%\marginpar{File: FFM-PRS}
\section{Polynomial recurrence sequences}
\label{sec:prs}

One type of MC-finite sequences is those obtained from polynomial recurrence sequences (PRS). These are sequences of vectors $\mathbf{a}_n \in \mathbb{Z}^k$ (for fixed $k$) where $\mathbf{a}_1,\dots,\mathbf{a}_c$ are given, and for $n > c$, each entry of $\mathbf{a}_n$ is a fixed polynomial with integer coefficients which depends on the $c$ preceding vectors $\mathbf{a}_{n-1},\dots,\mathbf{a}_{n-c}$. We call $k$ the \emph{dimension} and $c$ the \emph{depth}.

We can extract an integer sequence by looking at the first entry of all vectors. We call a sequence obtained in this way a polynomial recurrence sequence (PRS). These are MC-finite.

\begin{theorem} \label{thm:prs}
If $(a_n)_{n \in \mathbb{N}}$ is given by a PRS, then it is  MC-finite and $a_n \bmod m$ is FPT. Moreover, we can take $T(m) = m^{O(1)}$, and so $a_n \bmod m$ can be computed in polynomial time if $m$ is encoded in unary. 
%\textcolor{red}{(check!)}
\end{theorem}
\begin{proof}
Suppose that the PRS is obtained from the recurrence $(\mathbf{a}_n)_{n \in \mathbb{N}}$ of dimension $k$ and depth $c$.

There is an integer polynomial mapping $F\colon (\mathbb{Z}^k)^c \to (\mathbb{Z}^k)^c$ such that for all $n > c$,
\[
 F(\mathbf{a}_{n-c},\dots,\mathbf{a}_{n-1}) = (\mathbf{a}_{n-c+1},\dots,\mathbf{a}_n).
\]

Let $\mathbf{a}^{(m)}_n$ be obtained by taking every entry modulo $m$. We can project $F$ to a mapping $F^{(m)}\colon (\mathbb{Z}_m^k)^c \to (\mathbb{Z}_m^k)^c$ which can be computed explicitly. The sequence
\begin{align*}
&(\mathbf{a}^{(m)}_1,\dots,\mathbf{a}^{(m)}_n), \\
&F^{(m)}(\mathbf{a}^{(m)}_1,\dots,\mathbf{a}^{(m)}_n), \\
&F^{(m)}(F^{(m)}(\mathbf{a}^{(m)}_1,\dots,\mathbf{a}^{(m)}_n)), \\
\ldots
\end{align*}
is eventually periodic since it results from applying the mapping $F^{(m)}$ repeatedly, and $F^{(m)}$ has a finite domain. In fact, it must repeat within the first $m^{kc}+1$ steps by the pigeonhole principle.

We can now explicitly extract the representation of $(\mathbf{a}^{(m)}_n)_{n \in \mathbf{N}}$ 
which is
eventually periodic and which allows us to compute $a_n \bmod m$ in time polynomial in $n$.
\end{proof}

%--------------------------------------
%\newpage
%\marginpar{File: FFM-Specker}
\section{Specker's method}
\subsection{Density function}
Let $\tau$ be a finite relational vocabulary $\bar{R}= (R_1^{\rho(1)}, R_2^{\rho(2)}, \ldots, R_s^{\rho(s)})$
where $\rho(i)$ is the arity of $R_i^{\rho(i)}$ and $\rho = \max_i(\rho(i))$.
Let $\cC$ be a class of $\tau$-structures and $\cC(n)$ be the set of structures of the form $([n],\bar{R})$.
We define the {\em density function $d_{\cC}(n)$  of $\cC$} to be the number of ways the set $[n]$ can be made into a structure $([n], \bar{R}) \in \cC(n)$.

We want to study cases of $\cC$ where $d_{\cC}(n)$  is  (effectively) MC-finite.

\subsection{Structures of bounded degree}
Let $\fA = \langle A, \bar{R} \rangle$ be a $\tau$-structure. 
We define a symmetric relation $E_A$ on $A$, 
and call $\langle A , E_A \rangle$ the {\em Gaifman-graph of $\fA$.}
\begin{itemize}
\item
Let $a, b\in A$. 
$(a,b) \in E_A$ iff
there exists
a relation $R\in\bar{R}$ and some $\bar{a}\in R$ 
such that both $a$ and $b$ appear in $\bar{a}$ 
(possibly with other elements of $A$ as well).
\item
For any element $a\in A$, the {\em G-degree} of $a$
is the number of elements \\ $b\neq a$ for which  $(a,b) \in E_A$.
\item
$\fA$ is of {\em bounded G-degree $d$} 
if every $a \in A$ has G-degree at most $d$.
\item
We say that $\fA$ is {\em G-connected} 
if its Gaifman-graph is connected. 
\item
For a class of structures $\mathcal{P}$ we say it is of bounded degree $d$
(resp. connected)
iff all its structures are of bounded degree $d$ (resp. connected).
\end{itemize}
We note that the Gaifman graph may have loops.

\subsection{DU-matrices}
For two $\tau$-structures $\fA, \fB$ we denote by $\fA \sqcup \fB$ their disjoint union.
Let $\fA_i, i \in \N$ be an enumeration of all finite $\tau$-structures with  universe $[n]$ for some $n \in \N$ and let $\cC$ be a class of
$\tau$-structures. The DU-matrix of $\cC$, $DU(\cC)$ is the infinite  $(0,1)$-matrix 
where colums and rows are labeled $\fA_i$. We write $DU(\cC)_{i,j}$ for the entry $(\fA_i, \fA_j)$.
$$
DU(\cC)_{i,j} = \begin{cases} 1 & \fA_i \sqcup \fA_j \in \cC \\ 0 & \fA_i \sqcup \fA_j \not \in \cC \end{cases}.
$$
Let $rk(DU(\cC))$ be the rank of $DU(\cC)$ over the field $\Z/\Z_2$.

%--------------------------------------------------------
Below we show that there are uncountably many  graph properties of fixed DU-rank
Hence graph properties of fixed finite DU-rank are not necessarily definable even in Second Order Logic, or any other logical formalism
in which there are only countably many formulas.
\begin{proposition}
\label{ex:DU-rank}
Let $\cC$ be a class of finite $\tau$-structures. 
\begin{enumerate}[(i)]
\item
Assume all $\tau$-structures in $\cC$ are G-connected.
If we allow the empty set as a $\tau$-structure, then
$rk(DU(\cC)) =2$, but otherwise $rk(DU(\cC)) =0$.
%\marginpar{or
%otherwise $rk(DU(\cC)) =0$}
\item
Let $A \subseteq \N^+$ a set of non-negative integers and $\mathrm{Circ}$ the set of circles $C_n: n \in \N+$.
Define $\mathrm{Circ}_A = \{ C_i: i \in A\}$. Then  $rk(DU(\mathrm{Circ)} = rk(DU(\mathrm{Circ}_A)$ irrespective of the choice of $A$.
We conclude that there are uncountably many graph classes of the same DU-rank as for $\mathrm{Circ}$.
\item
Now let $\mathrm{Circ}_A^m$ be the class of graphs which consist of $k$ disjoint unions of graphs in $\mathrm{Circ}_A$.
If we allow the empty graph the rank $rk(DU(\mathrm{Circ}_A^m)$ is $m+1$ irrespective of $A$.
We conclude that there are uncountably many graph classes of the same DU-rank as for $\mathrm{Circ}^m$.
\end{enumerate}
\end{proposition}

From now on we allow structures to have an empty universe.
We next look at variations of  graph classes $\cC$ closed under disjoint union.
A graph property $\cC$  is a \emph{Gessel class of graphs} iff $\cC$ is 
closed under disjoint unions and closed under taking G-connected components.
\begin{facts}
\label{ex:gessel}
Let $\cC$ be a graph property.
\begin{enumerate}[(i)]
\item
If $\cC$ is a Gessel class of graphs  then $rk(DU(\cC))  \leq 2$.
\item
Let $G_1, G_2,  H$ be graphs and $H$ be connected. Then $G_1 \sqcup G_2$  does not contain $H$  as a subgraph (induced subgraph, minor)
iff neither $G_1$ nor $G_2$ contains $H$ as a subgraph (induced subgraph, minor).
\item
Let $\cC$ be a class of connected graphs and $\mathrm{Forb}_{sub}(\cC)$ 
the class of graphs which do not contain a graph $H \in \cC$ as a subgraph.
Then $\mathrm{Forb}_{sub}(\cC)$ 
is a Gessel class of graphs.
\item
The same holds if we replace $\mathrm{Forb}_{sub}(\cC)$ by 
$\mathrm{Forb}_{ind}(\cC)$, respectively $\mathrm{Forb}_{min}(\cC)$
where forbidden subgraphs are replace by forbidden induced subgraphs, respectively forbidden minors.
\item
Let $\cC$ a class of $\tau$-structures.
Assume that for all finite  $\tau$-structures $\fA, \fB$ we have
$\fA \sqcup \fB \in \cC$ iff both $\fA$ and $\fB$ are in $\cC$.
Then $rk(DU(\cC))$ is finite. 
%\\
%Note that this assumption does not imply that $\cC$ is closed under G-connected components.
%\end{enumerate}
%\end{facts}
\item
Let $\cD$ be a class of finite $\tau$-structures.
we denote by
$cl_{DU}(\cD)$ the closure of $\cD$ under disjoint unions.
If $\cC = cl_{DU}(\cD)$, then $rk(DU(\cC))$ is finite. 
\item
If $\cC$ consists of graphs which are the disjoint union of two equal-sized cliques, $rk(DU(\cC))$ is infinite.
%\item
%If $\cC$ consists of graphs on $n$ vertices with $n \equiv a \pmod{m}$  and which are the disjoint union of two cliques, 
%$rk(DU(\cC))$ is finite. In fact it is $m+1$.
\end{enumerate}
\end{facts}

%-----------------------------------------------------------------
Now we can state the main theorem on DU-rank:
\begin{theorem}[\cite{fischer2003specker,FiMa-ETH-2003}] 
\label{th:FiMa}
Let $\cC$ be a class of $\tau$-structures such that
\begin{enumerate}[(a)]
\item
every $\fA \in \cC$ is of G-degree at most $d$,  and
\item
$rk(DU(\cC)) =r$ is finite. 
\end{enumerate}
Then
\begin{enumerate}[(i)]
\item
the density function $d_{\cC}(n)$ of $\cC$ is MC-finite.
\item
If additionally all structures in $\cC$ are G-connected, then $d_{\cC}(n)$ modulo $m$ is ultimately $0$.
\end{enumerate}
\end{theorem}
In \cite{fischer2003specker} and its expanded version \cite{FiMa-ETH-2003}, 
the authors also exhibit a modular recurrence relation which depends only on the G-degree $d$ of $\cC$, 
the modulus $m$, and a fixed finite submatrix of $DU(\cC)$ containing $r$ rows. 
In the sequel we derive this recurrence relation in detail.

\subsection{Describing $DU(\cC)$}
In the sequel we assume that $\cC$ has finite DU-rank.
Two structures $\fA_i$ and $\fA_j$ are {\em DU-equivalent for $\cC$, equivalent for short}, if the columns $i$ and $j$ 
are identical in $DU(\cC)$. In this case we  write $\fA_i \sim \fA_j$.

\begin{observation}
\label{obs:equiv}
\begin{enumerate}[(i)]
\item
If $\fA_i \sim \fA_j$, then for every $\fB$ we have
$\fB  \sqcup \fA_i \sim \fB  \sqcup \fA_j $.
\item
As the rank of $DU(\cC)$ is finite, there are only finitely many equivalence classes.
\end{enumerate}
\end{observation}

Let $\cC_d$ be the class of structures of $\cC$ whose maximal degree is $\leq d$.
Let $\cD_0, \cD_1, \ldots, \cD_s$ be an enumeration of all DU-equivalence classes of $\cC_d$ with $\cD_0$ the only equivalence class which contains a structure
with maximal degree $\geq d+1$.
%Furthermore let $\cC_d = \{\fA \in \cC: 
%$$
%\cD_0 = \{ \fA : \forall \fB
%$$
Let $\cE = (\cD_1, \ldots, \cD_s)$ and let $d_{\cE}(n)$ be the vector-valued function
$d_{\cE}(n) = (d_{\cD_1}(n), \ldots, d_{\cD_s}(n))$ .

\begin{observation}
Let $DU(\cC)$ of finite DU-rank.
Modulo $m$ there are only finitely many possible values for $d_{\cE}(n)$.
\end{observation}

We shall use the following finite set of parameters of the matrix $DU(\cC)$:
\begin{enumerate}[(i)]
\item
The bound $d$ on the  G-degree of the structures in $\cC$.
\item
The rank $DU(\cC)$.
\item
Representatives $\fB_j$ for each of the equivalence classes $\cD_j$.
%With these we also get an algorithm for checking whether  $\fA \in \cD_j$.
\item
For fixed $m \in \N$ the finitely many values of 
$d_{\cE}(n) = (d_{\cD_1}(n), \ldots, d_{\cD_s}(n))$.
\end{enumerate}

In the case of, say, Gessel classes of graphs, these parameters are easily available.
To compute $d_{\cE}(n)$ in general we also need an algorithm for checking whether a finite $\tau$-structure $\fA$ belongs to $\cC$.
A  further algorithm for checking whether two $\tau$-structures are DU-equivalent for $\cC$ can then be obtained 
from the parameters and the algorithm above. 

\subsection{Orbits}

%\begin{defi}
Given a permutation group $G$ that acts on a set $A$
(and in the natural manner acts on models over the universe $A$),
the {\em orbit} in $G$ of a model ${\mathfrak A}$ with the universe
$A$ is the set
${\rm Orb}_G({\mathfrak A})=\{\sigma({\mathfrak A}):\sigma\in G\}$.
%\end{defi}

For $A'\subset A$ we denote by $S_{A'}$ the group
of all permutations for which $\sigma(u)=u$ for every $u\not\in A'$.
The following lemma is useful for showing linear congruences
modulo $m$.

\begin{lemma}
\label{dvlm}
Given ${\mathfrak A}$,
if a vertex $v\in A-A'$ has exactly $d$ neighbors in $A'$, then
$|{\rm Orb}_{S_{A'}}({\mathfrak A})|$ is divisible by ${|A'|\choose d}$.
\end{lemma}

\begin{proof}
Let $N$ be the set of all neighbors of $v$ which are in $A'$,
and let $G\subset S_{A'}$ be the subgroup
$\{\sigma_1\sigma_2:\sigma_1\in S_N\wedge\sigma_2\in S_{A'-N}\}$; in other
words, $G$ is the subgroup of the permutations in $S_{A'}$
that in addition send all members of $N$ to members of $N$.
It is not hard to see that
$|{\rm Orb}_{S_{A'}}({\mathfrak A})|
 ={|A'|\choose|N|}|{\rm Orb}_G({\mathfrak A})|$.
\end{proof}

The following simple observation is used to enable us to require in advance
that all structure in $\mathcal C$ have a degree bounded by $d$.

\begin{observation}\label{obs:bd}
We denote by ${\mathcal C}_d$ the class of all
members of $\mathcal C$ that in addition have bounded degree $d$.
If $\mathcal C$ has a finite $DU$-index, then so does ${\mathcal C}_d$.
%\noproof
\end{observation}

\subsection{Other parameters of the recurrence relation}
In order to show that $d^{(d)}_{\mathcal C}(n)$ is ultimately periodic
modulo $m$, we show a linear recurrence relation modulo $m$
on the vector function $(d_{\mathcal E}(n))_{\mathcal E}$ where
${\mathcal E}$ ranges over all other equivalence classes with respect
to ${\mathcal C}_d$.

The following parameters will be used in the recurrence relation:
\begin{enumerate}[(i)]
\item
$c = m\cdot(d!)$ which  depends only on the modulus $m$ and the G-degree $d$.
We note that for every $t\in{\mathbb N}$ and $0<d'\leq d$,
$m$ divides ${tc \choose d'}$. 
\item
$t= t(n)=\lfloor\frac{n-1}c\rfloor$.
\end{enumerate}

This, together with Lemma \ref{dvlm}, allows
us to prove the following.
\begin{lemma}\label{rclm}
Let ${\mathcal D}\neq{\mathcal N}_\phi$ be an equivalence class
for $\phi$, that includes the requirement of the maximum
degree not being larger than $d$. Then
$$
d_{\mathcal D}(n) \equiv
 \sum_{\mathcal E}a_{{\mathcal D},{\mathcal E},m,(n \bmod C)}
   d_{\mathcal E}(C\lfloor\frac{n-1}C\rfloor) \pmod{m},
$$
for some fixed appropriate $a_{{\mathcal D},{\mathcal E},m,(n \bmod C)}$.
\end{lemma}

\begin{proof}
Let $t=\lfloor\frac{n-1}C\rfloor$.
We look at the set of structures in ${\mathcal D}$ with the universe $[n]$,
and look at their orbits with respect to $S_{[tC]}$.
If a model ${\mathfrak A}$ has a vertex $v\in[n]-[tC]$
with neighbors in $[tC]$, let us denote the number of its neighbors
by $d'$. Clearly $0<d'\leq d$, and by Lemma \ref{dvlm} the size
of ${\rm Orb}_{S_{[tC]}}({\mathfrak A})$ is divisible by
${tC \choose d'}$, and therefore it is divisible by $m$. Therefore,
$d_{\mathcal D}(n)$ is equivalent modulo $m$ to the number of structures
in ${\mathcal D}$ with the universe $[n]$ that in addition have
no vertices in $[n]-[tC]$ with neighbors in $[tC]$.

We now note that any such structure can be uniquely
%structure ${\mathfrak A}$ with the universe $[n]$
%and no vertices in $[n]-[tC]$ with neighbors in $[tC]$ can be uniquely
%written as ${\mathfrak B}(1|{\mathfrak C})$ where ${\mathfrak B}$ is
%any structure with the universe $[n-tC+1]$ in which the vertex $v=1$
%has no neighbors, and ${\mathfrak C}$ is any structure over
written as ${\mathfrak B}\sqcup{\mathfrak C}$ where ${\mathfrak B}$ is
any structure with the universe $[n-tC]$,
and ${\mathfrak C}$ is any structure over
the universe $[tC]$. We also note using Observation \ref{obs:equiv}(i) that
the question as to whether ${\mathfrak A}$ is in ${\mathcal D}$ depends
only on the equivalence class of ${\mathfrak C}$ and on
${\mathfrak B}$ (whose universe size is bounded by the constant $C$).
By summing over all possible ${\mathfrak B}$ we get the required
linear recurrence relation (cases where
${\mathfrak C}\in{\mathcal N}^{(d)}_{\mathcal C}$ do not enter this sum
because that would necessarily imply
${\mathfrak A}\in{\mathcal N}^{(d)}_{\mathcal C}\neq{\mathcal D}$).
%${\mathfrak A}\not\in{\mathcal N}^{(d)}_{\mathcal C}$).
\end{proof}

\subsection{Proof of the main thoerem}

\begin{proof}[Proof of Theorem \ref{th:FiMa}:]
(i)
We use
Lemma \ref{rclm}:
%implies Theorem \ref{bdth}: 
Since there is only a finite
number of possible values modulo $m$ to
the finite dimensional vector $(d_{\mathcal E}(n))_{\mathcal E}$,
the linear recurrence relation in Lemma \ref{rclm} implies ultimate
periodicity for $n$'s which are multiples of $C$. From this
the ultimate periodicity for other values of $n$ follows,
since the value of $(d_{\mathcal E}(n))_{\mathcal E}$
for an $n$ which is not a multiple of $C$ is linearly related
modulo $m$ to the value at the nearest multiple of $C$.

(ii)
Finally, if 
all structures are connected we use Lemma
\ref{dvlm}.
Given ${\mathfrak A}$, connectedness implies that
%\marginpar{EF: Corrected}
there exists a vertex $v\in A'$ that has neighbors in $A-A'$. Denoting
the number of such neighbors by $d_v$, we note that
$|{\rm Orb}_{S_A'}({\mathfrak A})|$ is divisible by ${|A'|\choose d_v}$,
and since $1\leq d_v\leq d$ (using $|A'|=tC$) it is also divisible by $m$.
%every vertex $v\in A-A'$ has 
%neighbors in $A'$, say exactly $d_v$  many.
%Then
%$|{\rm Orb}_{S_A'}({\mathfrak A})|$ is divisible by all the values
%${|A'|\choose {d_v}}$.
%By our choice $|A'| = [tC]$, the above is divisible by $m$,
This makes the total number of models divisible by $m$
(remember that the set of all models with $A=[n]$ is a disjoint
union of such orbits), so
$f^{(d)}_{\mathcal C}(n)$ ultimately vanishes modulo $m$.
\end{proof}

%\subsection{Unbounded G-degree}

%--------------------------------------
%\newpage
%\marginpar{File: FFM-Specker-unbounded}
\subsection{Unbounded G-degree}
In the case of density functions for classes $\cC$ of unbounded G-degree the situtation is quite different.
In 2023, \cite{fischer2024extensions}, the following is shown:

\begin{theorem}
\label{th:ternary}
There is a class $\cD$ of $\tau$-structures where $\tau$ consists of one ternary relation such that
\begin{enumerate}[(i)]
\item
$\cD$ has unbounded G-degree,
\item
$DU(\cD)$ has finite DU-rank, and
\item
$d_{\cD}(n)$ is not MC-finite.
\end{enumerate}
\end{theorem}
Already in 2003, E. Fischer \cite{ar:Fischer02} constructed a class of $\tau'$-structures $\cD'$ 
with the same properties as $\cD$
where $\tau$ consists of one quaternary relation.

\begin{remark}
The original Theorem \ref{th:ternary} actually proved that $\cD$ is definable in First Order Logic $\FOL$.
Already in \cite{pr:BlatterSpecker84} it was noticed that $\FOL$-definability implies finite DU-rank.
One way to see this is by using the Feferman-Vaught Theorem for $\FOL$, 
see \cite[Theorem 14.4]{fischer2011application}.
\end{remark}

However, in the original paper \cite{pr:BlatterSpecker84}, C. Blatter and E. Specker did not use the DU-rank but
the substitution rank $Subst(\cC)$ of classes of $\tau$-structures $\cC$ where $\tau$ only contains finitely many relations
of arity at most $2$. They called it the {\em index} of $\cC$.
Their main result is:
\begin{theorem}[C. Blatter and E. Specker, 1981]
\label{th:SB}
Let $\cC$ be a class of $\tau$-structures where $\tau$ only contains finitely many relations
of arity at most $2$. Assume that $\cC$ has finite substitution rank.
Then $d_{\cC}(n)$ is MC-finite.
\end{theorem}

For $\tau$-structures $\fA$ where $\tau$ only contains finitely many relations
we define
\begin{enumerate}[(i)]
\item
A pointed $\tau$-structure $(\fA,a_0)$ is a $\tau$-structure with a distinguished element $a_0$.
We also allow $a_0 \not \in A$.
\item
Let $\fB$ another $\tau$-structure.
$Subst((\fA,a_0),\fB))$ is the structure obtained by replacing $a_0$ by $\fB$ and connecting all elements $a \in A$
with $R(a, a_0)$ for some binary $R \in \tau$ with all the elements $b \in B$.
\item
If $a_0  \not \in A$, the structure $Subst((\fA,a_0),\fB))$ is just the disjoint union of $\fA$ and $\fB$.
\item
The matrix $Subst(\cC)$ has columns labeled by $\tau$-structures and rows labeled by pointed $\tau$-structures.
The entry for $(\fA,a_0),\fB)$ is $1$ if  $Subst((\fA,a_0),\fB)) \in \cC$ and $0$ otherwise.
\item
The substitution rank $rk(Subst(\cC)$ is the rank of $Subst(\cC)$ over $\Z_2$.
\end{enumerate}

\begin{observations}
\begin{enumerate}[(i)]
\item
There is no obvious definition of $Subst((\fA,a_0),\fB)$ for ternary relations which works to obtain an analogue of Theorem
\ref{th:SB}.
\item
If $\tau$ only contains finitely many relations
of arity at most $2$, we have that
$rk(DU(\cC)) \leq rk(Subst(\cC))$, because $DU(\cC)$ is a submatrix of $Subst(\cC)$.
\end{enumerate}
\end{observations}

The proof of Theorem \ref{th:SB} also uses orbits but is much more complicated. A transparent  proof may be found in E. Specker's later
exposition in \cite{specker1990application} or in \cite[Section 16]{fischer2011application}.

%--------------------------------------
%\newpage
%\marginpar{File: FFM-ternary}

A {\em ternary structure} is given by the pair $\mathcal{A}=(A, B)$ where $A$ is a set and $B \subseteq A^3$.
A {\em ternary property} is a class of ternary structures closed under isomorphism.
Theorem \ref{th:FiMa} holds for ternary properties of bounded G-degree.
However, the combinatorial part of the Specker-Blatter theorem Theorem \ref{th:SB} does not hold for ternary properties 
of unbounded G-degree due to Theorem \ref{th:ternary}.
One immediate reason for this is the fact that the substitution rank for ternary  properties
is not well defined. However, this does not preclude the possibility that another general condition 
on classes of ternary structures $\cT$
exists in order to establish MC-finiteness.

\subsection{Examples of ternary structures and properties}

Natural finite ternary structures occur when a binary operation is represented by a ternary relation. Examples are groups and semi-groups.
The Gaifman graph of a binary operation as a ternary operation is a complete graph, therefore such a ternary property  is of unbounded G-degree,
provided it contains arbitrarily large ternary structures. The number of labeled ternary structures on $[n]$ is bounded by $2^{n^3}$.
In the literature the number of finite groups is usually counted up to isomorphisms. 

\begin{problem}
What is the number of ternary structures on $[n]$ which are finite groups?
Is it MC-finite?
\end{problem}

The most studied classes of ternary structures are motiviated by the study of the {\em betweenness relation} in elementary geometry.
There are many axiomatizations of betweenness relations having various geometric applications in mind.
For an encyclopedic survey the reader may consult \cite{pambuccian2011axiomatics}.
In the case where we study {\em finite} ternary structures the betweenness relations usually comes from partial orders, or more specifically,
from lattices.  In \cite{pitcher1942transitivities} the betweenness relations derived from lattices are extensively studied.

Another class of ternary structures are the {\em ternary spaces}
studied in \cite{hedlikova1983ternary}.
A ternary structure $\mathcal{A}$ is a {\em ternary space} if it satisfies the following axioms $\mathcal{T}$.
\begin{description}
\item[$T_1$:] 
$\forall a,b,c (B(a,b,c) \rightarrow B(c,b,a))$
\item[$T_2$:]
$\forall a,b,c (B(a,b,c) \wedge B(a,c,b)  \rightarrow b=c)$
\item[$T_3$:]
$\forall a,b,c, d (B(a,b,c) \wedge B(a,c,d)  \rightarrow B(b,c,d))$
\item[$T_4$:]
$\forall a,b,c,d (B(a,b,c) \wedge B(a,c,d)  \rightarrow B(a,b,d))$
\end{description}

\begin{problem}
What is the number of labeled ternary spaces?  Is it MC-finite?
\end{problem}

\subsection{Showing MC-finitness for ternary properties}

What do we have to know about a ternary property  $\mathcal{P}$ in order to show that  $d_{\mathcal{P}}(n)$ is MC-finite?

Here are some applicable methods:
\begin{description}
\item[Bounded degree:]
In \cite{fischer2003specker} it is shown that if all the structures in $\cC$ are of G-degree at most $d$,
then definability of $\cC$ in $\CMSOL$ implies MC-finiteness.
\item[PRS:]
In \cite{filmus2023mc} it is shown that if $d_{\cC}(n)$ satisfies a polynomial recurrence sequence, then it is MC-finite.
It may be difficult to find such a recurrence.
The method was briefly discussed in Section \ref{sec:prs}.
\item[Reducing to binary:]
In some cases there is a bijection of 
the ternary relations on $[n]$  of a ternary structure in  $\cT$ with the binary relations on $[n]$
of structures on $[n]$ in a class of structures $\cC$. In such a case (and similar cases), the density function $d_{\cT}(n)$ is the same as $d_{\cC}(n)$
where $DU(\cC)$ is of finite DU-rank.
This is the case for the betweenness relation of certain classes of partial orders.
\end{description}

\begin{problem}
Formulate and prove an analogue of Theorem \ref{th:SB} for ternary properties.
\end{problem}

%--------------------------------------
%\newpage
%\marginpar{File: FFM-conclu}
\section{Conclusions}

Theorem \ref{th:SB} uses the finite rank of $Subst(\cC)$ to prove MC-finiteness of $d_{\cC}(n)$.
In \cite{pr:BlatterSpecker84,specker1990application} it is shown that if $\cC$ is definable in  Monadic Second Order Logic $\MSOL$, then 
$\cC$ has finite substitution rank. This still holds if $\cC$ is definable in $\CMSOL$, the extension of $\MSOL$ by 
modular counting quantifiers, \cite{fischer2003specker}.
The proof of this uses the Feferman-Vaught Theorem for $\MSOL$, respectively for $\CMSOL$.

In its simplest form for disjoint unions  and for $\FOL$ it says:
\begin{theorem}[Feferman-Vaught Theorem for Disjoint Unions]
\label{th:fm}
%There is no elementary function $g$ such that the following holds for all finite $tau$-structures $\fA, \fB, \fC$: 
If $\fA$ and $\fB$ satisfy the same  $\tau$-sentences of $\FOL$ of quantifier rank at most $q$, then 
$\fA \sqcup \fC$ and $\fB \sqcup \fC$ 
satisfy the same first-order sentences of quantifier rank $q$.
\end{theorem}
The corresponding theorems for $\MSOL$ and $\CMSOL$ also hold, see \cite{makowsky2004algorithmic}.

From this one can derive that the DU-rank is finite.
However, such a proof gives very bad upper bounds for the DU-rank.
If $\cC$ is defined by a $\tau$-formula of $\MSOL$ of quantifier rank $q$,
the bound is roughly speaking determined by the number of non-equivalent  $\tau$-formulas of $\MSOL$ of quantifier rank $q$.

Another formulation of Theorem \ref{th:fm} is as follows:
\begin{theorem}
Let $\fA, \fB$ be two $\tau$-structures.
Assume $\phi$ is an $\FOL(\tau)$ sentence of quantifier rank $q$. Then there are $\FOL$ sentences 
$\psi_i: i \in I$, 
each of quantifier rank $q$,
and there is a boolean function $B$ in variables $x_i, y_i: i \in I$ such that the following holds:
$$
\fA \sqcup \fB \models \phi \text{ iff  }
B(a_i, b_i: i \in I) =1
$$
where $a_i, b_i$ are the truth value of $\fA \models \psi_i$ respectively $\fB \models \psi_i$. 
\end{theorem}
We call $\psi_i: i \in I$ a reduction sequence of $\phi$ and $|I|=I(\phi)$ its {\em FV-width}.
In this case the bound obtained for the DU-rank of $DU(\phi)$ is related to the  FV-width of $\phi$.

In fact, in \cite[Theorem 6.1]{dawar2007model} it is suggested that no better bounds can be obtained using the Feferman-Vaught Theorem.
They show a related result for the length of the formulas rather than the quantifier rank.

\begin{theorem}
\label{th:dawaretal}
There is no elementary function $g$ such that the following holds for all trees $\fA, \fB, \fC$: 
If $\fA$ and $\fB$ satisfy the same first-order sentences of length at most $g(\lambda)$, then 
$\fA \sqcup \fC$ and $\fB \sqcup \fC$ 
satisfy the same first-order sentences of length at most $\lambda$.
\end{theorem}

\begin{problem}
Can one prove a theorem similar to Theorem \ref{th:dawaretal} which estimates the minimal FV-width of $\phi$?
\end{problem}

They also note that if all the structures of $\cC$ are of G-degree bounded by $d$, the situation is much  better.
In this case there is an elementary function $g'$ and $g'$ is bounded by an iterated exponential function of at most $5$ iterations.

In this paper we discussed effective MC-finiteness in purely combinatorial terms.
In Examples \ref{ex:DU-rank} we have seen that in many cases the DU-rank of a class $\cC$ is very low and can be determined easily,
although using the Feferman-Vaught Theorem would only give a very large upper bound for it.
Similarly, one can construct examples with very low substitution rank.
For example the class of finite connected graphs is definable in $\MSOL$ by a formula of quantifier rank $3$,
hence has finite substitution rank. But it is easily seen that its substitution rank is $2$.
We leave it to the reader to find more examples.

\begin{problem}
Find more combinatorial %(non-logical) 
properties of $\cC$ which allow to determine the DU-rank or the substitution rank efficiently.
\end{problem}

%--------------------------------------
\subsection*{Acknowledgements}
The proof of Theorem \ref{th:FiMa} was taken verbatim from \cite{fischer2003specker,FiMa-ETH-2003}. It also appears in \cite{fischer2011application}.
\bibliography{FFM}
\bibliographystyle{alpha}
\end{document}